\documentclass{article}
\usepackage[utf8]{inputenc}
\usepackage{graphicx}
\usepackage{amsmath}
\usepackage[colorlinks=true,citecolor=red,linkcolor=blue,%
urlcolor=RubineRed,pdfpagetransition=Blinds,pdftoolbar=true,pdfmenubar=false]{hyperref}
\usepackage{color}

\usepackage{amssymb}
\usepackage{amsthm}
\usepackage{amsfonts}
\usepackage{eucal}
\usepackage{xspace}
\newtheorem{theorem}{Theorem}

\newtheorem{lemma}[theorem]{Lemma}

\newtheorem{remark}[theorem]{Remark}

\newcommand{\R}{\mathbb{R}}
\newcommand{\Z}{\mathbb{Z}}
\newcommand{\C}{\mathbb{C}}
\newcommand{\N}{\mathbb{N}}

\begin{document}
\title{Numerical approximation of the potential in the two-dimesional inverse scattering problem}
\author{ \textsc{Juan A. Barceló$\,^\ast$}\and \textsc{Carlos Castro}\thanks{Departamento de Matem\'atica e Inform\'atica, ETSI Caminos, Canales y Puertos, Universidad Polit\'ecnica de Madrid, Campus Ciudad Universitaria,
C/ Profesor Aranguren 3, 28040 Madrid, Spain. E-mails: {\tt juanantonio.barcelo@upm.es,} {\tt carlos.castro@upm.es.}}
\and \textsc{Juan M. Reyes}\thanks{School of Computer Science $\&$ Informatics, Cardiff University, Queen's Buildings, 5 The Parade, Roath, Cardiff CF24 3AA, United Kingdom. E-mail: {\tt reyes.juanmanuel@gmail.com.}} }
 \date{}

\maketitle

\abstract{We present an iterative algorithm to compute numerical approximations of the potential for the Schr\"odinger operator from scattering data. Four different types of scattering data are used as follows: fixed energy, fixed incident angle, backscattering and full data. In the case of fixed energy, the algorithm coincides basically with the one recently introduced by Novikov in \cite{Nov}, where some estimates are obtained for large energy scattering data. The numerical results that we present here are consistent with these estimates.}

\section{Introduction}

We consider the following scattering problem appearing in quantum physics,
\begin{equation} \label{eq_1}
\left\{ \begin{array}{ll}
(-\Delta +V(x)-k^2) u(x)=0, \qquad x\in \R^2, \\
u=u_i+u_s
\end{array} \right.
\end{equation}
where $u_i=u_i(x,\theta,k)=e^{ikx \cdot \theta}$ is the incident wave, with wave number $k\in\R$ and direction of propagation $\theta\in S^{1}$. The potential $V(x)\in L^\infty(\R^2)$ is assumed to be compactly supported and $u_s(x,\theta,k)$ is the scattered wave that satisfies the Sommerfeld radiation condition at infinity,
\begin{equation} \label{eq_2}
\frac{\partial u_s}{\partial r} (x,\theta,k)-iku_s(x,\theta,k)=o(r^{-1/2}), \qquad \mbox{ as $r\to \infty$} ,
\end{equation}
with $r=|x|$. From this condition we can deduce the following asymptotics for $u_s$,
\begin{equation} \label{eq_3}
u_s(x,\theta,k) =\frac{e^{ikr} u_\infty (\frac{x}{r},\theta,k)}{(kr)^{1/2}}  + o(r^{-1/2}),
\end{equation}
where the function $u_\infty (\theta',\theta,k)$ depending on $k\in\R$, the incident angle $\theta\in S^{1}$ and the reflecting angle $\theta'\in S^{1}$, is known as the {\em scattering amplitude} or {\em far field pattern} and constitutes the data of the inverse problem (see Ch. 5 in \cite{Ek} for details).

We are interested in recovering the quantum mechanical potential $V(x)$ from the far field data $u_\infty (\theta',\theta,k)$ for some values of $\theta'\in S^1$, $\theta\in S^1$ and $k\in\R$.

Let us state the problem in an equivalent integral formulation. Under the above conditions on $V(x)$, the system (\ref{eq_1})-(\ref{eq_2}) is equivalent to the so-called Lippmann-Schwinger equation
\begin{equation} \label{eq_35}
u(x,\theta,k)=e^{i k x\cdot \theta} + \int_{\R^2}\Phi(k|x-y|)V(y)u(y,\theta,k)\; dy
\end{equation}
where $\Phi(r)=\frac{i}{4}H_0^{(1)}(r)$ and $H_0^{(1)}$ is the Hankel function of the first kind of order zero.
On the other hand, the potential and the far field pattern are related by the integral equation (see \cite{ER}, \cite{Ek}),
\begin{eqnarray}
&&u_\infty(\theta',\theta,k)= \int_{\R^2} e^{-ik\theta'\cdot x} \,V(x)\; u(x,\theta,k)\; dx. \label{eq_4b}
\end{eqnarray}
Therefore, the inverse scattering problem can be stated as follows: Given the far field data $u_\infty (\theta',\theta,k)$, for some values of $\theta'\in S^1$, $\theta\in S^1$ and $k\in\R$, find a compactly supported $V(x)\in L^\infty(\R^2)$ that satisfies (\ref{eq_4b}), where $u$ is the solution of (\ref{eq_35}).

Note that $u_\infty$ depends on three variables while $V$ only depends on two. Thus, in principle, $V(x)$ could be recovered from the partial knowledge  of $u_\infty(\theta',\theta,k)$ on a suitable 2-dimensional submanifold of $S^1\times S^1\times \R$. The usual choices are: fixed energy $k=k_0$, fixed incident angle $\theta=\theta_0$ and backscattering $\theta'=-\theta$.


A number of algorithms have been proposed to recover the potential from the scattering data. Most of them are based on the Neumann-Born series for $V$ described in \cite{JK} (see also \cite{P1} and  \cite{P2}) which roughly consists in substituting (\ref{eq_35}) into the right hand side of (\ref{eq_4b}) iteratively. At each step a new multilinear term in $V$ appears in the right hand side of (\ref{eq_4b}). Assuming that $V$ can be written in power series with respect to a small parameter $\varepsilon$, i.e. $V=\sum_{n\geq 0} \varepsilon^n \widetilde{V}^{n+1}$ and making equal the same powers of $\varepsilon$ one easily obtains an iterative formula for $\widetilde{V}^n$.  The first term $\widetilde{V}^1$ of the series provides a somehow linearization of the inverse problem and it is known as the Born approximation of the potential. Other strategies are based on perturbation methods (see \cite{DO})
or more direct inversion algorithms as those proposed by Novikov (\cite{N1}, \cite{N2}, \cite{N3}, \cite{N4} and \cite{NS}). These have been implemented numerically in \cite{BBM} and \cite{BAR}, (see also the references therein).

In this paper we investigate numerically an iterative algorithm to approximate the potential based on a suitable fixed point iteration on the integral formula that defines the Born approximation, which we denote by $V_B$. We obtain a sequence of approximations $\{V^n\}_{n=1}^\infty$ where the new approximation $V^{n+1}$ can be deduced from $V^n$ and its associated scattered field $u^n_s$ by solving a single inverse Fourier transform.  To validate the algorithm we consider a numerical approximation of the inversion formula for the Fourier transform.  A convergence result for this numerical method to approximate the inverse Fourier transform is given in Section \ref{section:convergence}, using a technique which may be understood as the Nyquist–Shannon sampling theorem applied to the Fourier transform. As applications we obtain, on one hand, estimates for the numerical approximation of the Born approximation and, on the other hand, numerical evidences that the sequence $\{V^n\}_{n=0}^\infty$ converges and it provides a good approximation of the potential $V(x)$ in few iterations.

For scattering data at fixed energy $k=k_0$, the algorithm is similar to the one introduced by Novikov in \cite{Nov}, where the convergence of the iterative process is investigated from the theoretical point of view. The main difference with the algorithm described here is that in \cite{Nov}, at each iteration, the approximation of the potential is modified by a low-pass filtering/cutting process that we briefly explain at the end of Section 2 below. The convergence result in \cite{Nov} is stated for any dimension $d\geq 2$ and for smooth potentials $V$ (more than $d$ derivatives). The numerical results described here suggest that the convergence result in \cite{Nov} could be true also for less regular potentials ($V\in L^\infty$) and without the filtering-cutting process.

The rest of the article is organized as follows. The iterative algorithm to recover the potential is given in Section 2. The
numerical method to approximate both the Born approximation and the subsequent approximations $V^n$ is described in Section 3. An analysis of the convergence for the numerical approximation of $V_B$ in terms of the mesh step is performed in Section 4. In Section 5 we show how the scattering data was simulated from a potential example, and numerical experiments for both the Born approximation and the sequence $V^n$. Finally, some conclusions are presented in Section 6.

\section{The iterative algorithm to approximate the potential}

In this section we present the iterative algorithm to approximate the potential. To introduce it we first describe the Born approximation in detail. When substituting the equality $u=u_i + u_s$ into the right hand side of (\ref{eq_4b}) we obtain,
\begin{eqnarray} \nonumber
u_\infty(\theta',\theta,k) &=& \int_{\R^2} e^{-ik(\theta'-\theta)\cdot x} \,V(x) dx\\&&+\int_{\R^2} e^{-ik\theta'\cdot x}\, V(x)\;u_s(x,\theta,k) \; dx. \label{eq_4}
\end{eqnarray}
The Born approximation $V_B$ to $V$ is defined formally as the solution of (\ref{eq_4}) when neglecting the last nonlinear term, i.e. the solution of the linear problem,
\begin{equation} \label{eq_5}
u_\infty(\theta',\theta,k)= \int_{\R^2} e^{-ik(\theta'-\theta)\cdot x} V_B(x) dx.
\end{equation}

However, this identity is not consistent. Given $\xi\in \R^2$, the right hand side is constant for those values $(\theta',\theta,k)$ in the set
\begin{equation} \label{eq_6}
G_\xi =\{ (\theta',\theta,k)\in S^1\times S^1 \times \R, \mbox{ such that }k(\theta'-\theta) = 2\pi \xi \} ,
\end{equation}
while the left hand side of (\ref{eq_5}) does not satisfy this compatibility condition necessarily. Therefore, a proper definition of the Born approximation $V_B$ requires also a strategy to select, for a given $\xi\in \R^2$, some specific values $(\theta'(\xi),\theta(\xi),k(\xi))\in G_\xi$ in a unique way.

Following the methodology in reference \cite{PSS}, in order to define $V_B$ one has to choose two open subsets $\mathcal{M},\Omega$ of $S^1\times S^1\times \R$ and $\R^2$, respectively, such that the operator $\phi :\mathcal{M}\mapsto\Omega$ defined by $\phi(\theta',\theta,k)=(1/2\pi) k(\theta'-\theta)$, is an isomorphism. This way, from the definition \eqref{eq_5}, it follows that $u_{\infty}(\phi^{-1}(\xi)) = (\mathcal{F}\, V_B)(\xi)$ for $\xi\in \Omega$, where $\mathcal{F}$ denotes the Fourier transform operator. Hence, provided that $\R^2\setminus \Omega$ has zero measure, one can easily obtain $V_B$ by inverting the Fourier transform,
\begin{eqnarray}
\nonumber V_B(x)&=&\int_{\Omega} e^{i2\pi x\cdot\xi} u_{\infty}(\phi^{-1}(\xi))\, d\xi\\ \label{eq_6.1}
 &=&\int_{\Omega} e^{i2\pi x \cdot \xi} u_\infty(\theta'(\xi),\theta(\xi),k(\xi))\; d\xi.
\end{eqnarray}

Given $\xi=(\xi_1,\xi_2)\in \Omega$, we determine the values $(\theta'(\xi),\theta(\xi),k(\xi))\in G_\xi$ for the usual types of inverse potential scattering:
\begin{enumerate}
\item Case 1: Fixed energy $V_B^k$. We fix $k>0$. Here, $\Omega = D(0,k/\pi)$, $\mathcal{M}=\{(\theta' , \theta, k)\, :\, \theta,\theta'\in S^1\}$. Define  $(\theta'(\xi), \theta(\xi),k)\in G_\xi$ as follows,
\begin{equation}
\left\{ \begin{array}{l}
\theta (\xi)=\Theta(\alpha)\frac{\xi}{|\xi|}, \mbox{ where }\Theta(\alpha)=\left(\begin{array}{cc} \cos\alpha & -\sin\alpha\\ \sin\alpha & \cos\alpha\\\end{array}\right),\\
\theta'(\xi)=\theta(\xi)+2\pi \frac{\xi}{k}.\\
\end{array} \right.
\end{equation}
Here, $\Theta(\alpha)$ performs the rotation of angle $\alpha(\xi)=\arccos(-\pi\frac{|\xi|}{k})$ with respect to the origin. Note that, in this case, we can only describe values of $\xi$ in the disk $|\xi|\leq k/\pi$. Indeed, the set $\R^2\setminus\Omega$ does not have zero measure. Thus, in order to define properly the Fourier transform of $V_B$, we extend the far field pattern $u_{\infty}(\theta',\theta,k)$ by zero outside $\mathcal{M}=\phi^{-1}(\Omega)$, which means that a low pass filter is applied to the Born approximation. In this case, $V_B$ can not be compactly supported.
\item Case 2: Fixed incident direction $V_B^\theta$. We fix $\theta\in S^1$, so $\Omega=\{\xi\in\R^2\, :\, \xi\cdot \theta\neq 0\}$ and $\mathcal{M}=\{(\theta',\theta,k)\, : \, \theta'\in S^1\setminus\{\theta\},\, k\in\R\setminus 0\}$. Define $k(\xi)$ and $\theta'(\xi)$ as follows,
\begin{equation}
\left\{ \begin{array}{l}
k (\xi)= -\pi |\xi|^2/(\xi\cdot \theta),  \\
\theta'(\xi)=\theta(\xi)+2\pi \frac{\xi}{k(\xi)} .
\end{array} \right.
\end{equation}
Note that, in this case, we can only describe values of $\xi$ outside the hyperplane $\xi\cdot  \theta=0$, but this is a zero measure set.
\item Case 3: backscattering  $V_B^{\theta=-\theta'}$. We fix $\theta'=-\theta$, so $\Omega=\R^2\setminus 0$, and choose $\mathcal{M} = \{(-\theta,\theta,k)\, :\, \theta\in S^1,\, k>0 \}$ (the case $k<0$ would be possible as well). Define $k(\xi)$ and $\theta(\xi)$ as
\begin{equation}
\left\{ \begin{array}{l}
k (\xi)=|\xi|\pi, \\
\theta(\xi)=-\pi \frac{\xi}{k(\xi)} = - \frac{\xi}{ |\xi|} .
\end{array} \right.
\end{equation}
\item Case 4: full scattering data $V_B^f$. In this case the Born approximation is defined as an average of the Born approximations for fixed incident angle, i.e.
\begin{equation} \label{eq_7}
V_B^f(x)=\frac1{2\pi}\int_{S^1}V_B^\theta(x)\; d\theta .
\end{equation}
\end{enumerate}

In general, all these Born approximations recover some properties of the potential $V(x)$ and, in particular, its singularities in the scale of weighted $L^p$, H\"{o}lder and Sobolev spaces, see \cite{PSS}, \cite{OPS}, \cite{R}, \cite{RV}, \cite{Re}, \cite{BFRV} and \cite{BFRV1}.

Once described the Born approximation to $V$ we introduce the iterative algorithm to improve this approximation, which is basically a fixed point iteration applied to (\ref{eq_4}). We consider the sequence of approximate potentials $\{ V^n\}_{n=0}^\infty$ defined recursively by,
\begin{equation} \label{eq_30}
\left\{ \begin{array}{l}
V^0=0, \\
\mbox{For $n=0,1,2,...$, $V^{n+1}$ is the solution of}\\
 u_\infty(\theta',\theta,k)= \int_{\R^2} e^{-ik(\theta'-\theta)\cdot x}\, V^{n+1}(x) dx\\
\quad +\int_{\R^2} e^{-ik\theta'\cdot x} \, V^n(x)\;u_s^n(x,\theta,k) \; dx,
\end{array} \right.
\end{equation}
where $u_s^n(x,\theta,k)$ is the scattered field associated to the potential $V^n$, i.e. $u_s^n=u^n-e^{ikx\cdot \theta}$ and $u^n$ is the solution of (\ref{eq_35}) with $V=V^n$. 
To simplify, we derive directly $u_s^n$ as the solution of the following Lippmann-Schwinger equation,
\begin{eqnarray} \nonumber 
u_s^n(x,\theta,k)&=&\int_{\R^2} \Phi(k|x-y|)V^n(y)e^{i k y\cdot \theta}\; dy  \\&&+ \int_{\R^2}\Phi(k|x-y|)V^n(y)u_s^n(y,\theta,k)\; dy, \label{eq_35b}
\end{eqnarray}
which is equivalent to (\ref{eq_35}).

Note that $V^1$  satisfies the equation (\ref{eq_6.1}) for the Born approximation (i.e. $ V^1=V_B$) whereas, for $n>1$, $V^n$ can be interpreted as an improved approximation of $V$, since it takes into account a better approximation of the second integral in the right hand side of (\ref{eq_4}). Similar iterative procedures have been previously used by A. Ruiz in \cite[Section 5]{R} for theoretical issues, where an alternative second approximation $V^2$ is constructed by plugging the Born approximation into the first nonlinear term of the Neumann-Born series, or by Y.M. Chen  and W.C. Chew in \cite{ChCh} where a similar idea is considered to improve the Born approximation in a related electromagnetic inverse scattering problem.

As for the Born approximation, the system (\ref{eq_30}) is not consistent in general and a suitable strategy is required to recover $V^{n+1}$ from (\ref{eq_30})-(\ref{eq_35b}). To be more precise, the new potential $V^{n+1}$ is obtained from the previous one $V^n$ through the Fourier transform $\mathcal{F}$,
\begin{eqnarray}
\nonumber\mathcal{F} \, V^{n+1}(\xi)&=& u_\infty(\theta'(\xi),\theta(\xi),k(\xi))\\&&- \int_{\R^2} e^{-ik\theta'(\xi)\cdot y} \, V^n(y)\;u_s^n(y,\theta(\xi),k(\xi)) \; dy, \label{eq_31}
\end{eqnarray}
where $(\theta'(\xi),\theta(\xi),k(\xi))=\phi^{-1}(\xi)$ is defined according to one of the strategies defined before (fixed $k$, fixed $\theta$, etc.), and $\xi\in\Omega$. In this way, we have a different sequence $\{V^n\}_{n\geq 0}$ for each one of the cases defined before.

The iterative algorithm introduced above is far from being justified from the theoretical point of view. 
A rigorous validation would require to address in particular the following issues:
\begin{enumerate}
\item[I1.] The existence of $V^{n+1}$.
\item[I2.] If I1 holds true, the convergence of the sequence $\{V^n\}_{n\geq 1}$.
\item[I3.] If I2 holds true, the identification of the limit $V^*$ with the original potential $V$.
\end{enumerate}
Concerning the first issue, the definition of $V^{n+1}$, with $n\geq 1$, requires previously to compute the solution $u_s^n$ to the Lippmann-Schwinger equation (\ref{eq_35b}) corresponding to $V^n$. We do not know whether $V^n$ satisfies the conditions guaranteeing the existence and uniqueness of solution to this equation. 

The second issue is also completely open. The third one is related to the uniqueness  for the potentials $V(x)$ satisfying (\ref{eq_4}) but only when $(\theta',\theta,k)$ satisfies the restrictions provided by the case we are considering. For instance, in the case of the backscattering ($\theta'=-\theta$), assuming that we can pass to the limit in (\ref{eq_30}), $V^*$ will satisfy
\begin{eqnarray}
&&u_\infty(\theta',\theta,k)= \int_{\R^2} e^{-ik\theta'\cdot x} V^*(x)\; u^*(x,\theta,k)\; dx,\label{eq_4bb}
\end{eqnarray}
when $\theta =-\theta'$ and $k>0$. Here $u^*$ is the solution of (\ref{eq_35}) with the potential $V^*$. Therefore $V(x)=V^*(x)$ if equation (\ref{eq_4bb}) determines $V$ uniquely. This is also a difficult question, although local and generic uniqueness have been proved  by Eskin and Ralston (\cite{ER1}).  The uniqueness problem for fixed energy has been solved for Bukhgeim \cite{B} (see also \cite{GN}) and for fixed incident direction is also an open problem (see \cite{S} and \cite{Sa}). If we use all scattering data, $u_\infty(\theta', \theta, k)$   determines   the potential uniquely (see \cite{Sa} and \cite{PS}).

As we said in the introduction, the algorithm presented here is closely related to the one introduced in \cite{Nov} for fixed energy scattering data. In this case, the sequence of approximations $\{W^n\}_{n\geq 0}$ starts also from $W^0=0$ and it is constructed iteratively as follows: we consider again formula (\ref{eq_31}) but replacing $V^n$ by $W^n$ and $V^{n+1}$ by $\tilde W^{n+1}$, an intermediate function. The new iteration $W^{n+1}$ is then obtained from $\tilde W^{n+1}$ in two steps: we first compute a suitable low-pass filter of $\tilde W^{n+1}$, and then multiply it by a  compactly supported cutoff function. With this strategy the issue I1 above is satisfied for $\{W^n\}_{n\geq 0}$, due to the fact that $W^n$ is compactly supported at each step (by construction), while the convergence of the sequence $\{W^n\}_{n\geq 0}$ is deduced from a careful analysis based on the filtering step. In this way, the following estimate is obtained,
\begin{equation} \label{eq_nov}
\| W^n - V\|_{L^\infty} \leq C_n k^{-\alpha_n},
\end{equation}
where $C_n$ is uniformly bounded, 
$$
\alpha_n=\left( 1-\left( \frac{r-d}{r}\right)^n\right)\frac{r-d}{2d},
$$
$d$ is the space dimension (in our case $d=2$), and $r>d$ is the required  number of derivatives in $L^1$ of the potential $V$. Note that this estimate does not imply the convergence of $W^n$ to $V$ as $n\to \infty$, for fixed $k$, but it provides a polynomial convergence rate as $k\to \infty$, for fixed $n$.  
We compare below our numerical estimate for $V^n-V$ with this one.

\section{Numerical approximation}

In this section we describe the numerical approach to find approximations of both the Born approximation defined by formula (\ref{eq_6.1}) and the sequence $V^n$ of iterative approximations defined by  (\ref{eq_31}). We first introduce some notation on the finite dimensional space and then we state the numerical versions of (\ref{eq_6.1}) and (\ref{eq_31}).

\subsection{Finite dimensional trigonometric space}
 Given $R>0$, we define
$$
G_R=\left\{ x=(x_1,x_2)\in \R^2 \; : \; |x_k|<R,\; k=1,2 \right\} .
$$
The family of exponentials
$$
\varphi_j(x)=e^{2i\pi\xi_j \cdot x}, \quad \xi_j=\frac{j}{2R}, \quad j=(j_1,j_2)\in\Z^2,
$$
constitutes an orthonormal basis on $L^2(G_R)$ with the norm
$$\| u \|_0^2=\frac1{(2R)^2}\int_{G_R}|u(x)|^2 dx.$$ We also introduce the space $H_\lambda=H_\lambda(G_R)$ which consists of $2R-$multiperiodic functions (distributions) having finite norm
$$
\| u \|_\lambda=\left( \sum_{j\in\Z^2} |{\underline \xi_j}|^{2\lambda} |\hat u(\xi_j)|^2 \right)^{1/2} ,
$$
with
$$
{\underline \xi_j}=\left\{ \begin{array}{ll} \xi_j, & (0,0)\neq j \in \Z^2 \\ 1, & j=(0,0), \end{array}  \right.
$$
and
$$
\hat u(\xi_j)=\int_{G_R} u(x)\overline{\varphi_j(x)} dx, \quad j\in \Z^2,
$$
the Fourier coefficients of $u$.

We now introduce a finite dimensional approximation of $H_\lambda$. Let us consider $h=2R/N$ with $N\in \N$ and a mesh on $G_R$ with grid points $jh$, $j\in \Z_h^2$ and
$$
\Z_h^2=\left\{ j=(j_1,j_2)\in Z^2 \; : \; -\frac{N}{2} \leq j_k < \frac{N}{2}, \; k=1,2 \right\}.
$$
We also consider $\mathcal{T}_h$ the finite dimensional subspace of trigonometric polynomials of the form
$$
v_h=\sum_{j\in Z_h^2} c_j \varphi_j, \qquad c_j\in \C.
$$
Any $v_h\in \mathcal{T}_h$ can be represented either through the Fourier coefficients
$$
v_h(x)=\sum_{j\in Z_h^2} \hat v_h(j) \; \varphi_j(x),
$$
or the nodal values
$$
v_h(x)=\sum_{j\in Z_h^2} v_h(jh) \; \varphi_{h,j}(x),
$$
where
$$
\varphi_{h,j}(x)=h^2\sum_{k\in \Z_h^2} e^{i\pi k \cdot (x-jh)/R}.
$$
For a given $v_h \in \mathcal{T}_h$, the nodal values $\bar v_h$ and the Fourier coefficients $\hat v_h$ are related by the discrete Fourier transform $\mathcal{F}_h$ as follows,
$$
\hat v_h =h^2\mathcal{F}_h \bar v_h, \qquad \bar v_h =\frac{1}{h^2}\mathcal{F}_h^{-1} \hat v_h,
$$
where, as usual, $\mathcal{F}_h$ relates the sequence $x(k)$ with $X(k)$ according to
$$
X(k)=\sum_{n=0}^{N-1} x(n)e^{-i2\pi n/N}, \qquad k=0,1,...,N-1 .
$$

The orthogonal projection from $H_\lambda$ to $\mathcal{T}_h$ is defined by the formula
$$
P_h v=\sum_{j\in Z_h^2} \hat v(\xi_j) \varphi_j,
$$
while the interpolation projection $Q_hv$ is defined, when $\lambda>1$, by
$$
Q_hv\in \mathcal{T}_h, \quad (Q_hv)(jh)=v(jh), \quad j\in \Z^2_h.
$$

\subsection{Finite dimensional setting}

Let $R>0$ such that supp$(V)\subset G_R$. We approximate the Born approximation in (\ref{eq_6.1}) by the following finite dimensional version: Find $V_{B,h} \in \mathcal{T}_h$ such that
\begin{equation} \label{eq_51}
 h^2 \mathcal{F}_h V_{B,h}(j)=u_\infty(\theta'(\xi_j),\theta(\xi_j),k(\xi_j)), \quad j\in \Z_h^2.
\end{equation}
Note that $V_{B,h}$ is computed from a single inverse discrete Fourier transform from the values of the far field $u_\infty(\theta'(\cdot),\theta(\cdot),k(\cdot))$ at the mesh points $\xi_j$ with $j\in \Z_h^2$. More precisely, the numerical approximation is obtained from the following process:

\bigskip

{\bf Algorithm 1:}
\begin{enumerate}
\item  Choose $h$ according to the mesh grid where we will compute the nodal values of $V_{B,h}$: $x_j=jh$, $j\in\Z_h^2$.
\item Construct the mesh $\xi_j=j/(2R)$ with $j\in\Z_h^2$.
\item Compute $u_\infty(\theta'(\cdot),\theta(\cdot),k(\cdot))$ at the mesh points $\xi_j$.
\item Invert the discrete Fourier transform to obtain the values of $V_{B,h}$ at the nodes $x_j$.
\end{enumerate}

Analogously, we compute the numerical approximations of $V^{n+1}$ in (\ref{eq_31}) by solving the following discrete formula for $V^{n+1}_h\in \mathcal{T}_h$
\begin{eqnarray}
\nonumber&&h^2 \mathcal{F}_h V^{n+1}_h(j)= u_\infty(\theta'(\xi_j),\theta(\xi_j),k(\xi_j))\\&&\quad - h^2 \sum_{k\in Z_h^2} e^{-i\theta'(\xi_j)\cdot kh} V^n_h(kh)\;u_s^n(\theta'(\xi_j),\theta(\xi_j),k(\xi_j)),  j\in \Z_h^2, \label{eq_42}
\end{eqnarray}
where $u_s^n(x,\theta(\xi_j),k(\xi_j))$ is obtained from a discrete formulation of the Lippmann-Schwinger equation (\ref{eq_35b}). In the experiments that we present in Section 5 below we follow the trigonometric collocation method introduced by G. Vainikko \cite{Va}. 

The iterative process is then as follows:

\bigskip

{\bf Algortihm 2:}
\begin{enumerate}
\item  Choose $h$ according to the mesh grid where we will compute the nodal values of $V_{B,h}$: $x_j=jh$, $j\in\Z_h^2$.
\item Construct the mesh $\xi_j=j/(2R)$ with $j\in\Z_h^2$.
\item Compute $V_{B,h}=:V^1_h$ following the Algorithm 1 above.
\item For $n=1,...,K$: $V_h^n \to V_h^{n+1}$
\begin{enumerate}
\item Compute $u_s^n(x,\theta(\xi_j),k(\xi_j))$ by approximating the Lippmann-Schwinger equation (\ref{eq_35b}) with $V=V^n_h$, for each mesh point $\xi_j$, $j\in \Z_h^2$.
\item Compute the right hand side of (\ref{eq_42}) for each mesh point $\xi_j$.
\item Invert the discrete Fourier transform to obtain the values of $V_{h}^{n+1}$ at the nodes $x_j$.
\end{enumerate}
\item End
\end{enumerate}

Note that each iteration requires $N^2$ performances of the  Lippmann-Schwinger equation solver, which is a costly computation. The good news is that, on a hand, these computations can be parallelized, and on the other hand, the results can be estimated by interpolation with the ones generated from a few computations on a coarse grid.

In the next section, we give estimates for the error of the
numerical approximation $V_{B,h}$ in (\ref{eq_51}) with respect to
a periodized version of the Born approximation $V_{B}$. As a
result, if $V_{B}$ has compact support in $G_R$, we have estimates
for the difference between the approximation $V_{B,h}$ and the
true $V_{B}$, since in this case the periodized version coincides
with the Born approximation $V_{B}$ in $G_R$. We do not know
whether the Born approximation preserves or not the compactness of
the potential's support, except for the fixed energy case, in
which a low pass filter is applied to $V_B$, and consequently
$V_B$ is not compactly supported.

\section{Convergence of the numerical approximation of the Born approximation}\label{section:convergence}

In this section, we prove a convergence result for the numerical
approximation of the Born approximation obtained by the Algorithm
1. The result is a direct consequence of a Lemma below that
provides an estimate for a discrete approximation of the inverse
Fourier transform.

\begin{theorem} \label{th_1}
Let $V_B$ the Born approximation of a potential $V(x)$ defined by
(\ref{eq_6.1}), according to one of the possible cases 1-4 defined
in the introduction (fixed energy, fixed incident angle, etc.). Consider a mesh of size $h$ as before in such a way that $\xi_j\in \Omega$ for all $j\in \Z_h^2$.
Let $V_B^\sharp$ the periodized version of $V_B$ defined as
\begin{equation}\label{eq_52e}
 V_B^\sharp(x)=\sum_{j\in \Z^2_h} V_B(x+2Rj).
\end{equation}
If $V_B^\sharp \in H^\lambda$ for some $\lambda>0$, then
\begin{equation} \label{eq_52f}
\| V_{B,h} -V_B^\sharp \|_0 \leq h^\lambda \| V_B^\sharp
\|_\lambda,
\end{equation}
where $V_{B,h}\in \mathcal{T}_h$ is the solution of (\ref{eq_51}).
\end{theorem}

\begin{remark}
The periodic version of $V_B$ given in (\ref{eq_52e}) will
coincide with $V_B$ in $G_R$ only if $V_B$ is compactly supported
in $G_R$. We do not know if this is the case, in general.
\end{remark}

The proof of Theorem \ref{th_1} is a direct consequence of the following lemma, which is an adapted version of the sampling theorem.

\begin{lemma} Let $q(x)\in L^2(\R^2)$ be such that
$$
\mathcal{F} q (\xi)=g(\xi), \quad \xi \in \R^2,
$$
where $g(\xi)$ is a continuous known function. Consider the finite dimensional approximation given by $q_h\in \mathcal{T}_h$, the solution of
\begin{equation} \label{eq_52bb}
h^2 \mathcal{F}_h q_h (j)=g(\xi_j), \quad j \in \Z_h^2.
\end{equation}
Let  $q^\sharp$ be the periodized version of $q$ given by
\begin{equation} \label{eq_52c}
 q^\sharp(x)=\sum_{j\in \Z^2_h} q(x+2Rj).
\end{equation}
In particular, if $q$ is compactly supported in $G_R$ then $
q^\sharp =q$ on $G_R$.

Then, if $q^\sharp \in H^\lambda$ for some $\lambda>0$,
\begin{equation} \label{eq_52b}
\| q_h - q^\sharp \|_0 \leq h^\lambda \|  q^\sharp \|_\lambda.
\end{equation}
\end{lemma}

{\bf Proof:} We divide the proof in several steps.

Step 1: Filtering. In this first step we define a low pass filter of $g(\xi)$, denoted by $g_{LP}(\xi)$, in such a way that
$$
supp \; (\mathcal{F}^{-1} g_{LP}) \subset G_R.
$$
We take
$$
g_{LP}(\xi)=\sum_{j\in \Z^2} g(\xi_{j})m(\xi-\xi_j),
$$
where 
$$
m(\xi)=\frac{\sin(2\pi R \xi^1)}{2\pi R\xi^1}\frac{\sin(2\pi R \xi^2)}{2\pi R\xi^2}, \qquad \mbox{ with }\xi=(\xi^1,\xi^2),
$$
and $\xi_j=j/(2R),$ $j\in\Z^2.$ 
Note that with this choice
$$
g_{LP}(\xi_j)=g(\xi_{j}), \qquad j\in\Z^2.
$$

Let us define $\tilde q \in L^2$ as the solution of
\begin{equation} \label{eq_40}
\mathcal{F} \tilde q (\xi)=g_{LP}(\xi), \quad \xi \in \R^2.
\end{equation}
By construction $\tilde q$ has compact support in $G_R$. Moreover, taken into account that the Fourier transform of $\chi_{G_R}(x)e^{i2\pi\xi_j\cdot x}$ is $(2R)^2m(\xi-\xi_j)$ we have,
\begin{equation} \label{eq_41bb}
\tilde q(x)=\mathcal{F}^{-1} g_{LP}(x)= \frac{1}{(2R)^2}\sum_{j\in \Z^2} g(\xi_{j}) e^{i2\pi\xi_j\cdot x} \chi_{G_R}(x),
\end{equation}
where $\chi_{G_R}(x)$ is the characteristic function of the two dimensional interval $G_R$.

\bigskip

Step 2: Cut and periodize. The main idea is to establish a
convergence result for $\tilde q$ in $G_R$ with the $H_\lambda$
metric. Therefore we define $ q^\sharp\in H_\lambda$ as the
$2R$-periodic extension of $\tilde q$ which coincides with $\tilde
q$ in $G_R$. We show that $P_h q^\sharp$ is in fact $q_h$ defined
by (\ref{eq_52bb}).

First of all, write the restriction of $\tilde q (x)$ to $G_R$ in Fourier series, i.e.
\begin{equation} \label{eq_41}
\tilde q (x)= q^\sharp (x)=\sum_{j\in \Z^2} c_j \varphi_j(x), \quad
x\in G_R.
\end{equation}
It is easy to see, from equation (\ref{eq_41bb}), that $c_j=g_{LP}(\xi_j)=g(\xi_j)$.

The projection $P_h  q^\sharp \in \mathcal{T}_h$ is given by
$$
P_h  q^\sharp(x)=\sum_{j\in \Z_h^2} c_j \varphi_j(x), \quad x\in \R^2.
$$
Then, we have for the nodal values of $P_h\tilde q$, $\overline{P_h\tilde q}$
$$
h^2 \mathcal{F}_h (\overline{P_h q^\sharp }) (j)=g(\xi_j), \qquad
j\in \Z_h^2,
$$
which is the formula that defines $q_h$. Therefore,
$$
q_h=P_h q^\sharp .
$$

\bigskip

Step 3. Now we estimate the norm in (\ref{eq_52b})
$$
\| q_h - q^\sharp \|_{L^2} = \| P_h q^\sharp -  q^\sharp \|_0 \leq
(4h)^{\lambda}  \|  q^\sharp\|_\lambda,
$$
where we have used that the projection operator $P_h$ satisfies, in general,
$$
\|P_hv-v\|_\lambda \leq \left(4h \right)^{\mu-\lambda} \| v\|_\mu, \qquad \mbox{for all $v\in H_{\mu}$, \; $\mu\geq \lambda$}
.
$$

\bigskip

Step 4. Finally, we prove that $ q^\sharp$  is given by  formula
(\ref{eq_52c}). As $\tilde q= q^\sharp $ in $G_R$ it is sufficient
to prove it for $\tilde q$.  Taking into account (\ref{eq_41bb})
we have,
\begin{eqnarray*}
\tilde q(x)&=&\mathcal{F}^{-1} g_{LP}(x)= \frac{1}{(2R)^2}\sum_{j\in \Z^2} g(\xi_{j}) e^{-i2\pi\xi_j\cdot x} \chi_{G_R}(x) \\&=&
\sum_{j\in \Z^2}q(x+2R j) \chi_{G_R}(x),
\end{eqnarray*}
where we have used the Poisson summation formula in the last
identity. This completes the proof of (\ref{eq_52c}).

\section{Numerical experiments}

In this section we show the efficiency of Algorithms 1-2 to approximate the potential, when inverting the Fourier transform with the numerical method described in the previous section. We divide this section into three subsections. The first one is devoted to the direct problem, i.e. we show how to construct synthetic far field data for a given potential in order to test our algorithms. This is also useful to apply the Algorithm 2.  The approximation of the potential from the far field data is discussed in the next two subsections, we consider separately the results with the Algorithm 1, that provides the Born approximation, and the Algorithm 2.

\subsection{Data simulation} \label{sec:data_sim}

We consider different situations according to the different cases
presented in the introduction. The far field data
$u_\infty(\theta'(\xi),\theta(\xi),k)$ is simulated by solving the
direct problem with the numerical code written by K. Knudsen, J.
L. Mueller and S. Siltanen to solve numerically the D-bar equation
for the two-dimensional Calder\'on problem. These codes are based
on the numerical approach to solve Lippmann-Schwinger equations
introduced by G. Vainikko \cite{Va}.

For each scattering type, consider the mesh $\{\xi_j = j/(2R)\, :\, j\in\Z^2_h\}$. Fix $j$ and write $\theta_j:= \theta(\xi_j)$, $k_j:= k(\xi_j)$. We compute a numerical approximation $\tilde{u}_s(x_l,\theta_j,k_j)$ to the solution $u_s(\cdot , \theta,k)$ of the problem (\ref{eq_1})-(\ref{eq_2}) for $\theta=\theta_j$, $k=k_j$, by adapting the aforementioned codes by S. Siltanen et al. to a $2R$-biperiodic version of the integral equation (\ref{eq_35}). To this end, we write (\ref{eq_35}) as follows
\begin{equation}\label{eq_43}
[I-\Phi_k\ast (V\cdot (\cdot))]u_s(\cdot ,\theta, k) = \Phi_k\ast (V e^{ik\theta\cdot (\cdot)}),
\end{equation}
where $\Phi_k(x):=\Phi (k|x|)$ denotes the fundamental solution to the Helmholtz equation, and periodize the equation (\ref{eq_43}) by replacing all its terms with their $2R$-biperiodic extension, cutting the Green's function smoothly and taking the convolution operator on the torus. We follow the spirit of Section 15.4 in the book \cite{MS} by J. L. Mueller and S. Siltanen. The fundamental solution $\Phi_k$ is computationally implemented through the Matlab function {\it besselh}.

The approximation $\tilde{u}_s(x_l , \theta_j , k_j)$ is computed on a grid $\{x_l = lh\, :\, l\in\Z^2_h\}$ of the square $G_R$, with $h=2R/N$. Assuming a potential supported in the unit disc $D(0,1)$, the condition $R>2$ is necessary for the periodization argument.

We have added a $5\%$ Gaussian noise to the synthetic scattering
data $\tilde{u}_s(x_l , \theta_j , k_j)$ to simulate possible
measurement errors and validate the robustness of the approach.
More precisely, we take $(1+0.05\, \mathcal{N}) \tilde{u}_s$ as
scattering data, where $\mathcal{N}\in \R^{N^2}$ is a sample of a
random vector whose elements are independent Gaussian random
variables with zero mean and unit standard deviation.

Next, the approximation of $u_{\infty}(\theta'(\xi_j) , \theta(\xi_j) , k(\xi_j))$ is generated from $\tilde{u}_s(x_l , \theta_j , k_j)$ via numerical quadrature using (\ref{eq_4}). Note that the values $\tilde{u}_s(x_l , \theta_j , k_j)$ are only reliable for $x_l\in D(0,1)$. Nevertheless, this is not a problem in (\ref{eq_4}), since $V$ is supported in $D(0,1)$.

We have chosen $R=2.1$ and as test potential $V=\chi_1+1.2\chi_2$ where $\chi_1$ is the characteristic function of the annulus $0.7<|(x_1,x_2)|<1$ and $\chi_2$ is the characteristic function of the square $|x_1|+|x_2|<0.3$ inside the annulus.

\subsection{The Born approximation: algorithm 1}

\begin{figure}
\centerline{ \includegraphics[width=7cm]{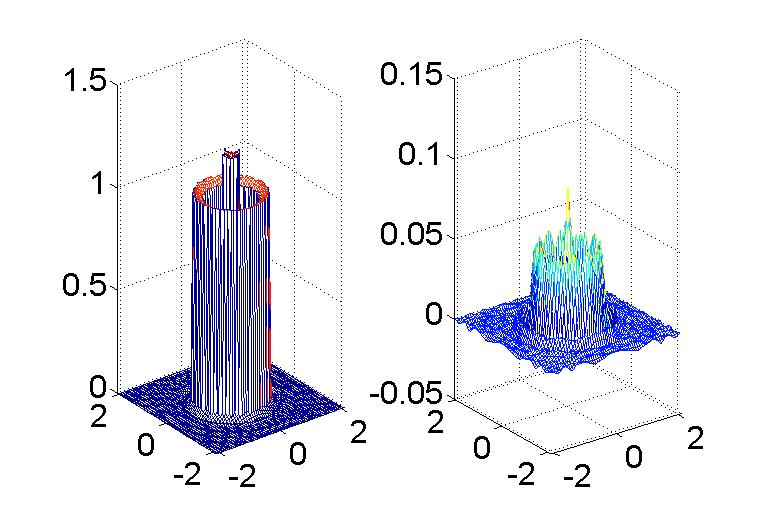} }
\caption{ Numerical approximation of the Born approximation for fixed energy $k=10$ and $N=2^6$. The left figure corresponds to the real part and the right one to the imaginary part. \label{fig_1}}
\end{figure}

\begin{figure}
\centering
\begin{picture}(330,130)
\put(116,175){real}
\put(181,175){imag}
\put(15,112){$V$}
\put(15,40){$V_{B,h}$ }
\put(34,-20){ \includegraphics[width=7cm]{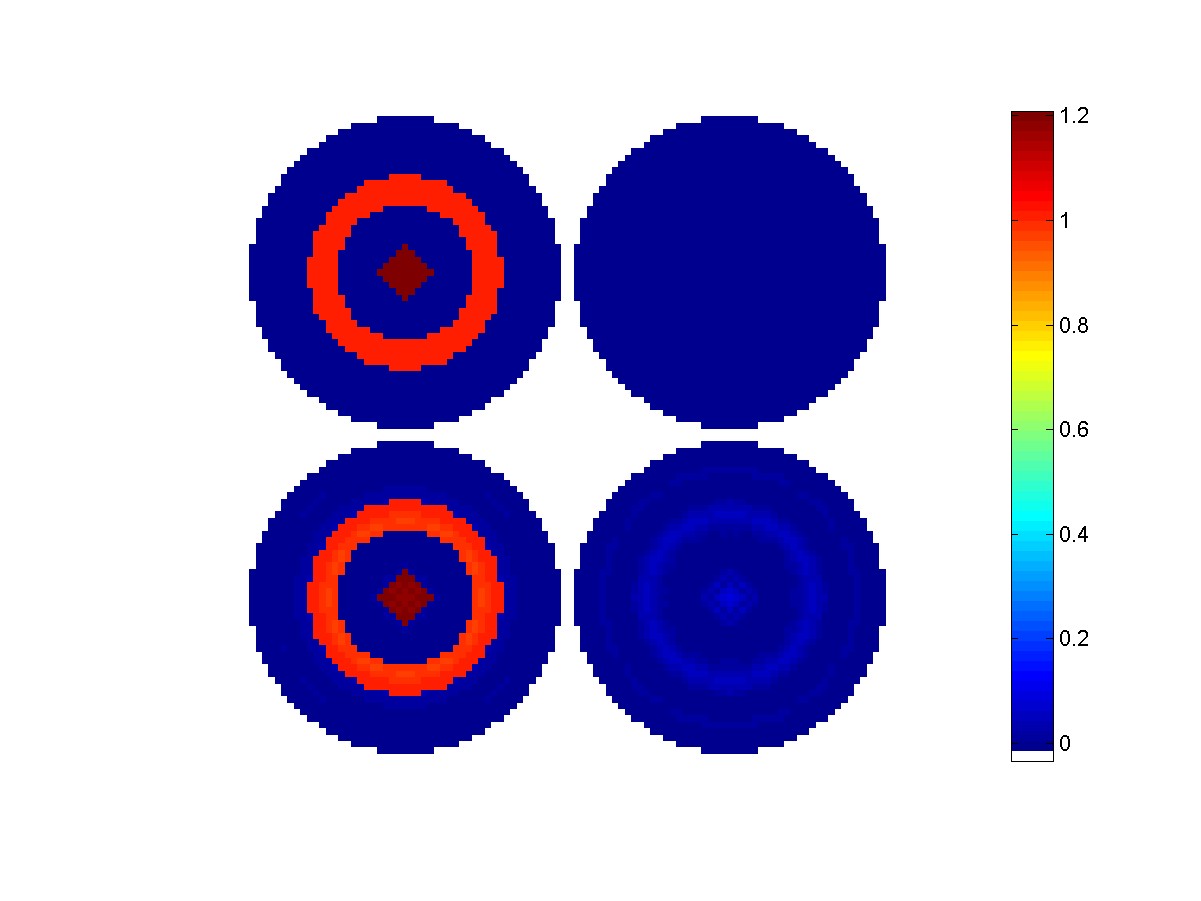} }
\end{picture}
\caption{\label{fig_1bis} Images of the true potential (top) and the numerical approximation of the Born approximation for fixed energy $k=10$ and $N=2^6$ (bottom). The left column corresponds to the real parts and the right column to the imaginary parts. For aesthetic reasons, only the values on the disc of radius 1.6 are shown. The colour scale is the same for the four pictures.}
\end{figure}
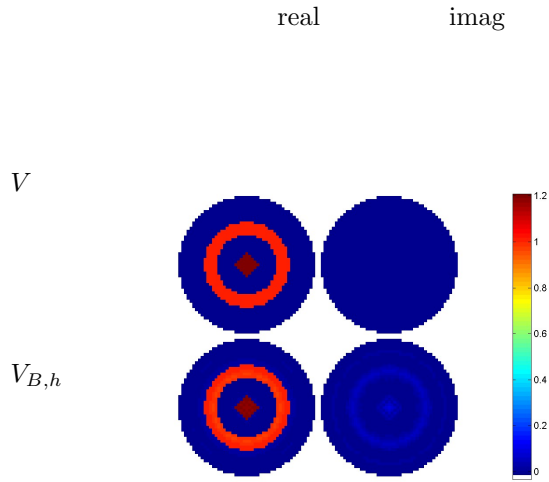

We first consider the fixed energy Born approximation. In this case, we can only compute $\theta(\xi)$ and $\theta'(\xi)$ for  $|\xi|\leq k/\pi$.  Combining the restriction $|\xi|<k/\pi$ and the fact that the mesh grid for $\xi$ corresponds to the interval $|\xi|<N/(4R)$ we deduce that we have only far field data for the whole meshgrid in (\ref{eq_51}) when
$$
 N<\frac{4Rk}{\pi}.
$$
If we choose $k=10$ for example, this bound is $N\sim 27$. Thus, the far field data can be computed in the whole mesh $\xi_j$ with $j\in\Z^2_h$ when $N\leq 27$. The right hand side in (\ref{eq_51}) is assumed to be zero for  $|\xi_j|\geq k/\pi$. In Figure \ref{fig_1} we show the real and imaginary parts of the reconstruction. Further, Figure \ref{fig_1bis} depicts the real and imaginary parts of this reconstruction, together with the same for the true potential, in order that the reader can visualize the accuracy of the approximation. Only the values on the disc of radius $1.6$ are shown. For comparison purposes, the same colour scale is used in all the pictures.  

\begin{figure}
\centerline{\includegraphics[width=7cm]{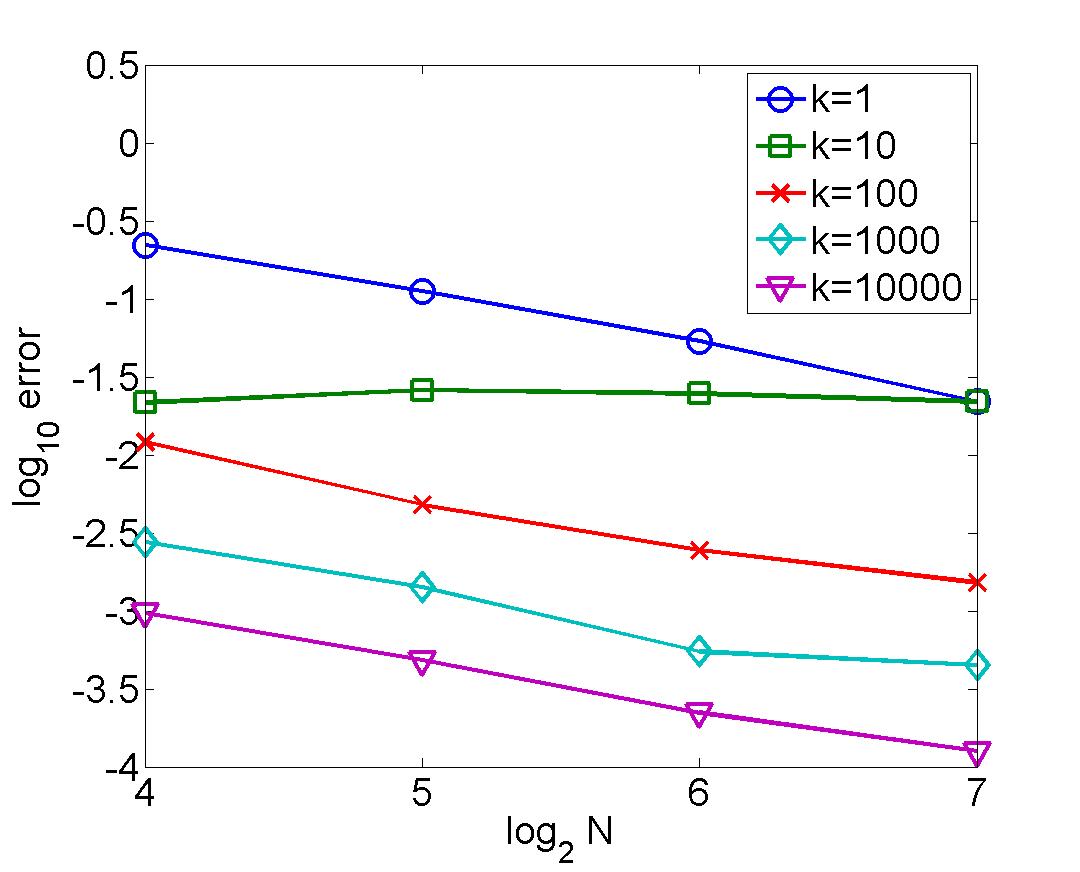}}
\caption{ Numerical approximation of the Born approximation for different values of $k$.  \label{fig_2}}
\end{figure}

Figure \ref{fig_2} shows the behavior of the error between the potential and the numerical approximation of the Born approximation in terms of $k$ and $N$. Here, and in the sequel, this error is computed with the following approximation of the $L^2$-norm,
$$
error^2=\frac{4R^2}{N^2}\sum_{j\in\Z^2_h}|V_{B,h}(x_j)-V(x_j)|^2.
$$
Note that this formula incorporates two types of errors: the one associated to the numerical approximation of $V_B$ and the distance of the Born approximation to the real potential $V$. There is no way to distinguish between these two errors since we do not know the Born approximation $V_B$ for this particular example.

From the results of Figure \ref{fig_2} we see that the error decreases for larger $k$ and larger $N$.

We now compare the rest of the Born approximations.

\begin{figure}
\begin{center}
\begin{tabular}{c}
 \includegraphics[width=7cm]{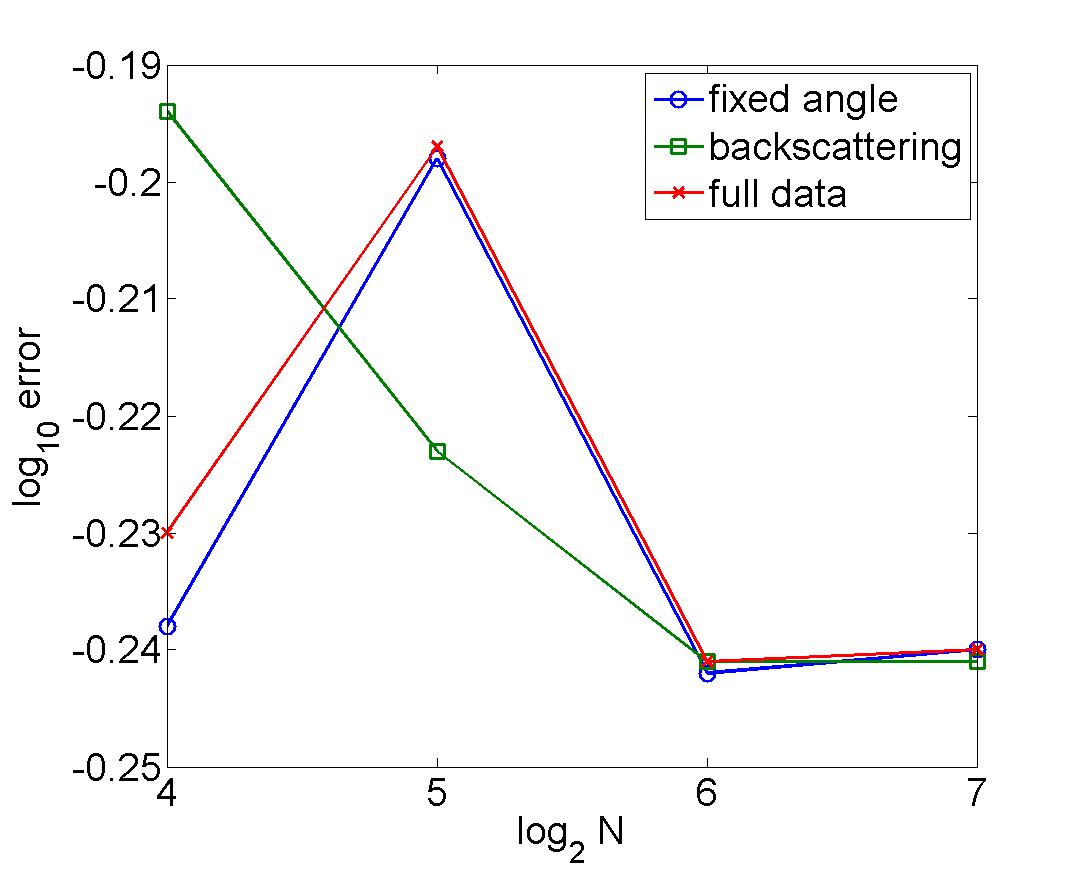}
\end{tabular}
\caption{ Distance between the potential and the real part of the different Born approximation: $\theta=\pi/4$ fixed, backscattering and full data cases, depicted with circles, squares and crosses, respectively. \label{fig_3}}
\end{center}
\end{figure}

In the case 2, fixed incident angle $\theta$, we can only compute $\xi$ outside the hyperplane $\xi\cdot \bar \theta=0$. The right hand side of (\ref{eq_51}) is then assumed to be zero when $\xi$ is in this hyperplane. The difference between the real potential and the Born approximation is very similar for any of the incident angles. Thus, we only show results for the case $\theta=\pi/4$ in Figure \ref{fig_3}.

In the case 4, full data, we have averaged the Born approximations
of 10 different incident angles equally distributed in $[0,2\pi)$. In
Figure \ref{fig_3} we compare the error between the potential and
the real part of these Born approximations for different mesh
sizes $N\times N = 2^M\times 2^M$, $M=4,5,6,7$. We see that, for
$M\geq 6$ the backscattering, fixed incident angle and full data generate
the same error roughly. However, the fixed-$k$ case analyzed
previously is significatively better than all these Born
approximations. We also observe that, in contrast with the
fixed-$k$ Born approximation, the error seems to stabilize after
$M=6$.

\begin{figure}
\begin{center}
\begin{tabular}{cc}
 \includegraphics[width=5cm]{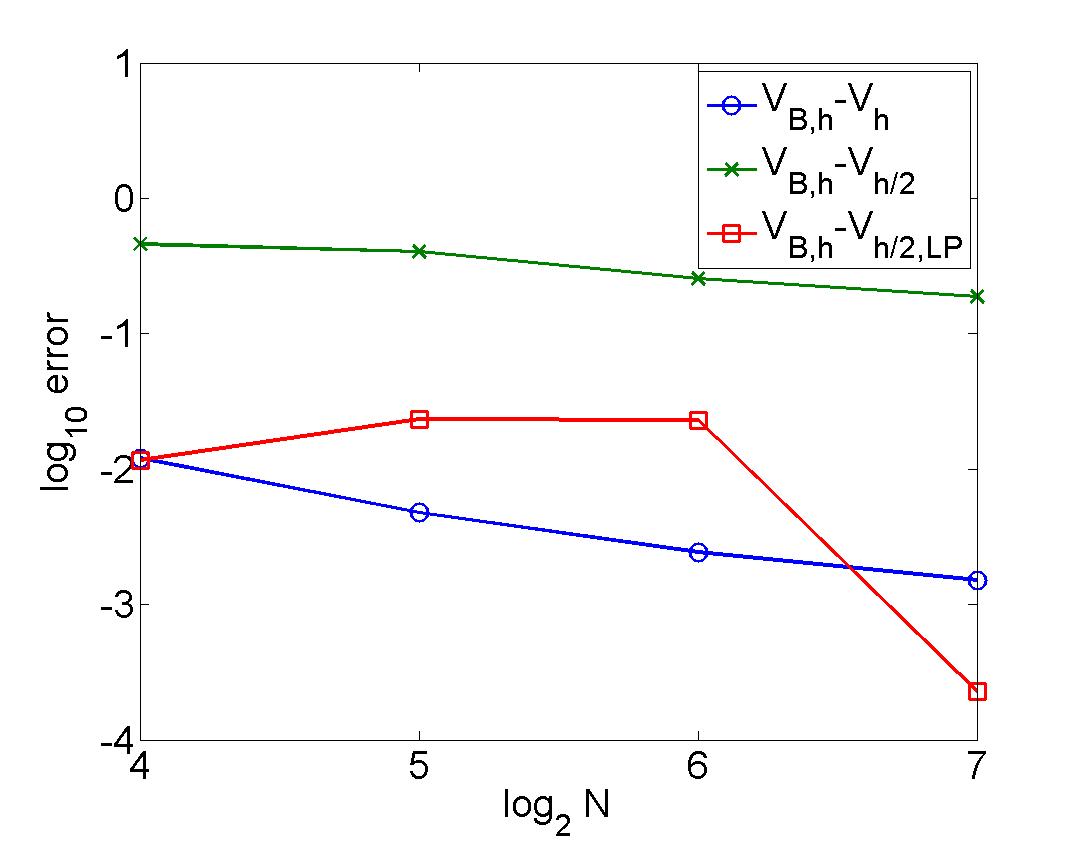} &
\includegraphics[width=5cm]{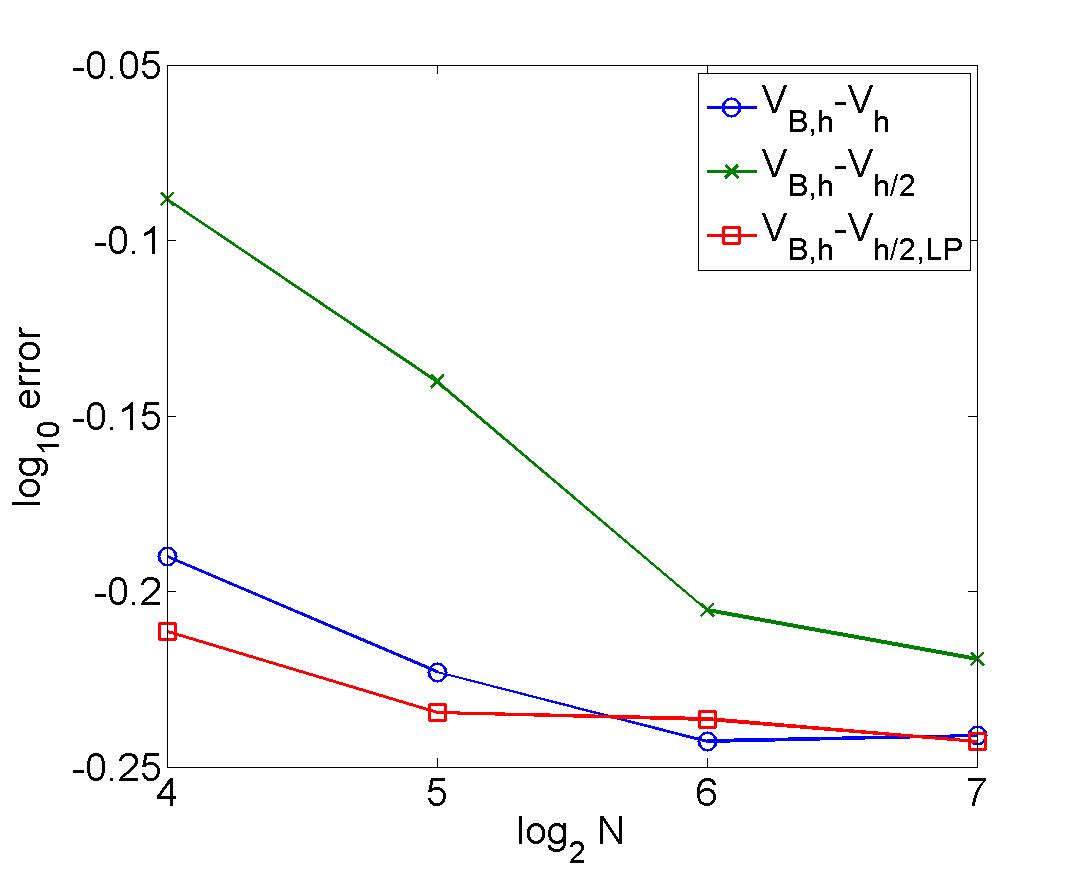}
\end{tabular}
\caption{ Distance between the numerical Born approximation of the potential and the potential when scattering data is computed in the same grid $V_{B,h}-V_h$, when a twice finer mesh is used to compute the scattering data $V_{B,h}-V_{h/2}$, and when it is compared with a low pass filter of $V$ in the finer mesh $V_{B,h}-V_{h/2,LP}$. The left figure corresponds to the case $k=100$-fixed while the backscattering case is in the right one.   \label{fig_alia}}
\end{center}
\end{figure}

In the simulations above the synthetic far field data are
computed in the same meshgrid as the one used to recover the potential. This is an important fact to avoid aliasing effects since we recover the potential from frequency data. If we use either a coarser or finer mesh to compute $u_\infty$ from the potential then it will contain necessarily an aliasing effect.  To illustrate this, in Figure \ref{fig_alia} we show the convergence of the numerical approximation for fixed $k=100$ and backscattering, as $N$ grows, when the data for the inverse problem (following Subsection \ref{sec:data_sim}) is computed in a mesh twice finer than the mesh used to recover $V_{B,h}$ with the Algorithm 1. We compare the error of the numerical Born approximation $V_{B,h}$ with respect to the real potential $V$ in the mesh of size $h$, and with respect to a low pass filter of the potential in the finer grid of size $h/2$ that removes half frequencies. As we see the latter error is smaller. 

This can be interpreted as the fact that we cannot recover anything better than a low pass filter of the potential $V$ which takes into account the frequencies represented in the frequency mesh given by $\xi_j$.

\subsection{The iterative algorithm}

In this section we show the efficiency of the Algorithm 2 to approximate the potential. Note that in this case, at each iteration, we have to simulate scattering data for the current approximate potential $V^n$ in order to recover the next approximation $V^{n+1}$ via the inverse Fourier transform. This is done following the idea explained in subsection \ref{sec:data_sim} above for simulating the synthetic data in the experiments. 

\begin{figure}
\begin{center}
\begin{tabular}{c}
 \includegraphics[width=7cm]{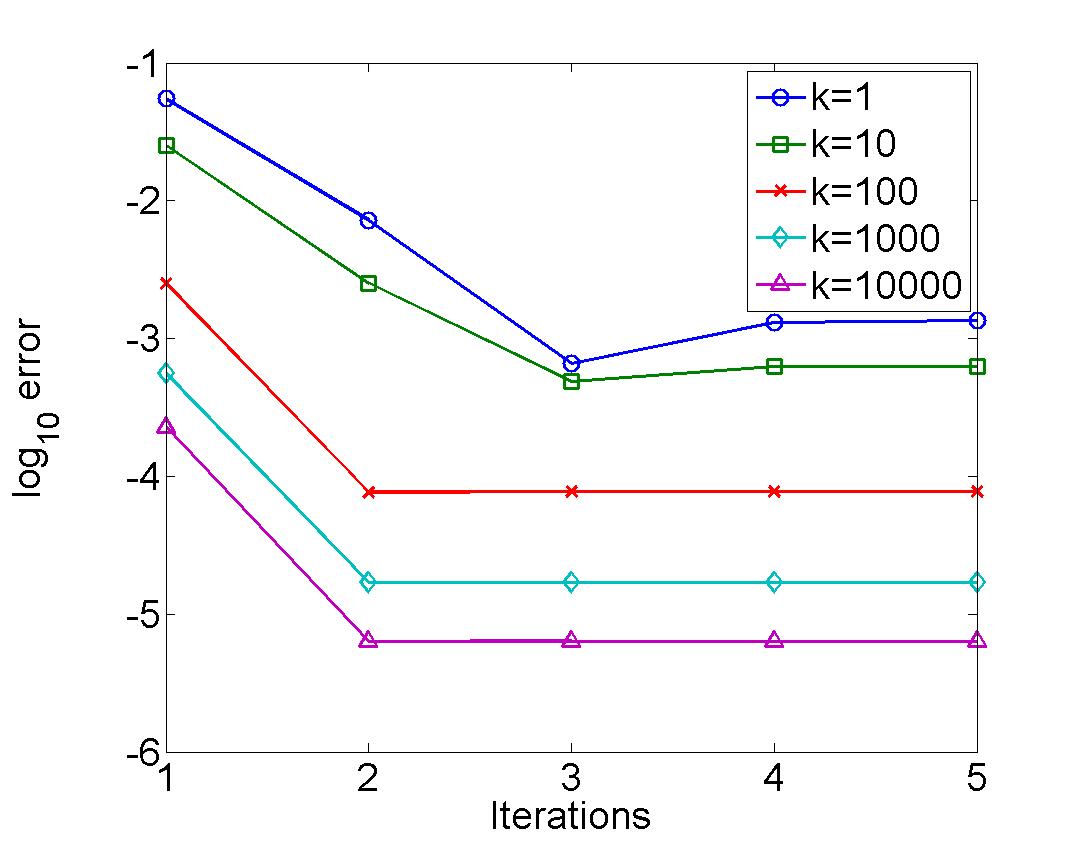}
\end{tabular}
\caption{ Distance between the potential and the real part of the numerical approximation when $k$ is fixed, for different values of $k$, in terms of the number of iterations. Here, the grid is taken with $N^2=2^{12}$ points. In the profiles, the circles, squares, crosses, rhombi and triangles refer to $k=1$, $k=10$, $k=100$, $k=1000$, $k=10000$, respectively,  as in Figure \ref{fig_2}. \label{fig_4}}
\end{center}
\end{figure}

In Figure \ref{fig_4} we show the difference between the real part of the approximation and the potential for different values of fixed $k$ and $N=2^6$ grid points in each variable. We see that the algorithm decreases in the first few iterations and then it stabilizes. Thus, for fixed $k>0$ we observe a convergence $V^n\to V^* \neq V$ but close to $V$.  

\begin{figure}
\begin{center}
\begin{tabular}{c}
 \includegraphics[width=7cm]{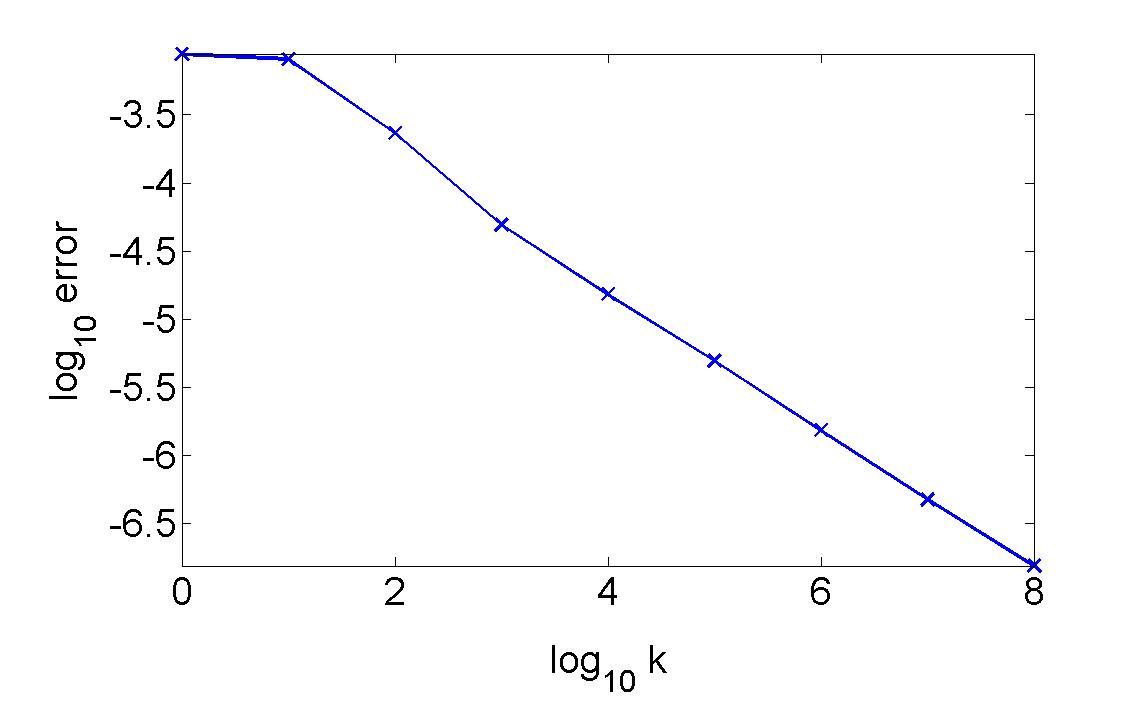}
\end{tabular}
\caption{ Error in $V^n-V$ versus $k$ after $n=6$ iterations. The grid has $N^2=2^{12}$ points. \label{fig_4b}}
\end{center}
\end{figure}

We also observe in Figure \ref{fig_4} that, as $k$ is larger, the approximation becomes better. The behavior is similar for other grid sizes. The convergence of the iterative method as $k\to \infty$ is also illustrated in Figure \ref{fig_4b} where we fix the number of iterations to $n=6$ and compare $V^6-V$ for different values of $k$. We obtain numerically the following convergence rate: 
\begin{equation} \label{eq rate_nov}
error \sim k^{-0.48}. 
\end{equation}
This rate is similar for other different initial data in $L^\infty$. The situation agrees with the estimate (\ref{eq_nov}), obtained by Novikov in \cite{Nov} for the modified sequence of approximated potentials $W^n$, in the sense that the convergence rate is a power of $k$. However, the exponent that we find numerically for $V^n-V$ is slightly better than the one given by estimate (\ref{eq_nov}). In fact, $\alpha_\infty=0.5$ when $r=4$ while (\ref{eq rate_nov}) holds for the $L^\infty$ potentials that we have tested.

\begin{figure}
\begin{center}
\begin{tabular}{c}
 \includegraphics[width=7cm]{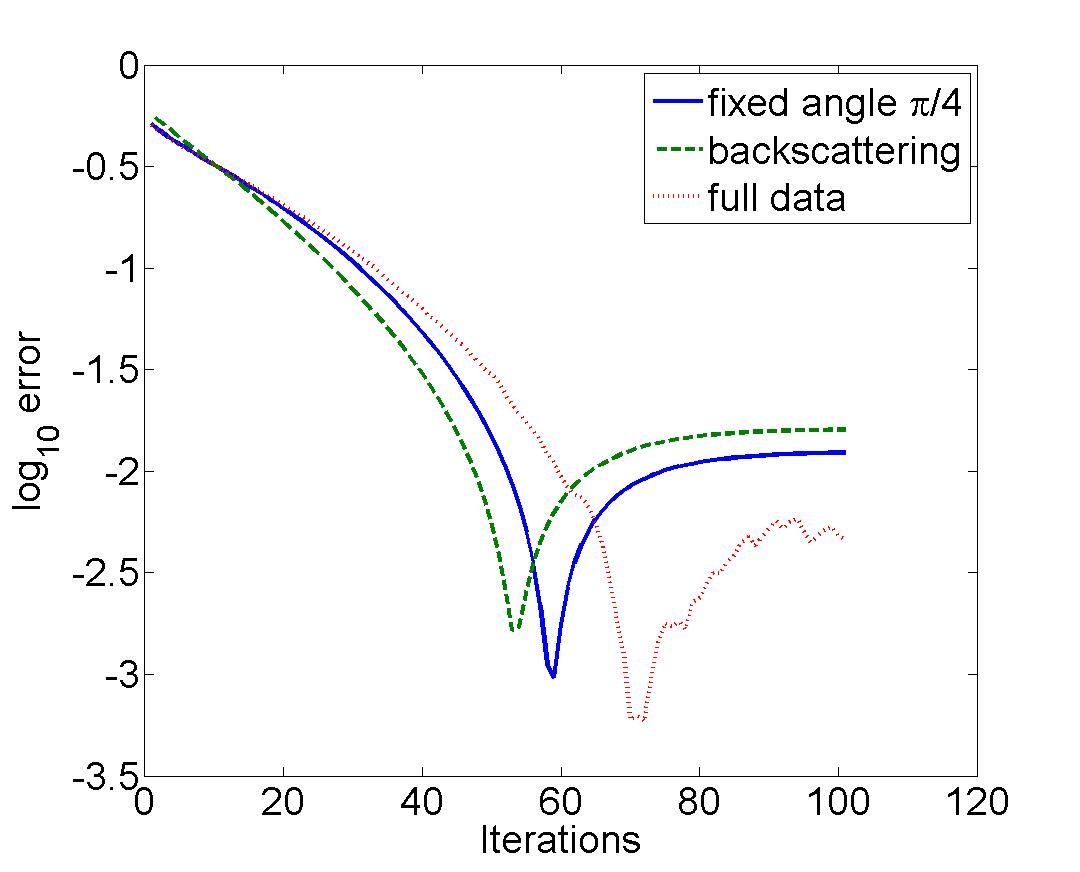}
\end{tabular}
\caption{ Distance between the potential and the real part of the numerical approximation in terms of the number of iterations, for fixed incident angle, backscattering and full data. An $N^2=2^{12}$-point mesh is considered. \label{fig_5}}
\end{center}
\end{figure}

In Figure \ref{fig_5} we compare the convergence for the cases 2-4
with $N=2^6$ grid points in each variable. We see that in the
three cases the distance to the original potential decreases until
a point from which it starts to increase. This is due to the fact
that we are considering noisy scattering data and the algorithm
converges to a ``noisy'' potential, which is not exactly the
original potential $V(x)$.

\section{Conclusions}

The inverse scattering problem for the Schr\"odinger equation consists of the recovery of the electrostatic potential $V$ from scattering data measurements modelled by the far field pattern $u_{\infty}$. The numerical approximation of this problem is investigated for different class of scattering data, namely, fixed energy, fixed incident angle, backscattering and full data.

Two different algorithms in 2D are considered. The Algorithm 1 performs an inversion of the discretized far field pattern using the discrete Fourier transform. In this way the so-called Born approximation is approximated. The Algorithm 2 goes beyond the Algorithm 1 through an iterative process solving successively the Lippmann-Schwinger (L-S) equation corresponding to the previous iteration result. The motivation of this iterative approach is the idea that these successive solutions to the L-S type equations approximate the scattered wave corresponding to the original potential, which is not proven and we propose as an open problem.

A convergence result is provided for the Algorithm 1. This task is especially difficult for the Algorithm 2 and is outside the scope of this work. The theoretical open problems that the Algorithm 2 pose, mentioned in the Introduction, are challenging and we hope they attract the attention of further articles.

Both algorithms are tested from noisy scattering data and $L^\infty$ potentials without rotational symmetry.

\section*{Acknowledgements}

The authors were supported by the project from Spain MTM2014-57769-C3-2-P (MINECO). The first and third authors were also supported by the project MTM 2011-02568 (MINECO). The second author has been supported by the project MTM2011-29306-C02-01 from the MICINN (Spain). J.~M. Reyes was additionally supported by the Academy of Finland (Finnish Centre of Excellence in Inverse Problems Research 2012–2017,
decision number 250215), the Engineering and Physical Sciences Research Council (EPSRC) from the United Kingdom, reference EP/K024078/1.

The authors thank Prof. Samuli Siltanen from the University of Helsinki for sharing with them his Matlab codes solving numerically Lippmann-Schwinger (L-S) type equations. J.~M. Reyes also thanks him for suggesting the idea of adapting his L-S equation solver used for the Calder\'on problem to the L-S equation of this work.

The authors thank the anonymous referee for pointing us the reference \cite{Nov}.

\end{document}